\documentclass[11pt]{article}
\usepackage{amsmath, amssymb, amsfonts, amsthm, latexsym} 
\usepackage{mathrsfs}
\usepackage{graphicx}
\usepackage{authblk}
\usepackage[usenames]{color}
\newtheorem{theorem}{Theorem}[section]
\newtheorem{lemma}[theorem]{Lemma}

\newtheorem{definition}[theorem]{Definition}

\newtheorem{corollary}[theorem]{Corollary}

\newtheorem*{notat*}{Notation}

\title{The Hamilton cycle space of random regular graphs and randomly perturbed graphs}

\author{Dan Hefetz \thanks{School of Computer Science, Ariel University, Ariel 40700, Israel. Email: {\tt danhe@ariel.ac.il}.}
\quad Michael Krivelevich \thanks{School of Mathematical Sciences, Tel Aviv University, Tel Aviv 6997801, Israel. Research supported in part
by NSF-BSF grant 2023688. Email: {\tt krivelev@tauex.tau.ac.il}.}}

\begin{document}
\maketitle
 
\begin{abstract}
The cycle space of a graph $G$, denoted $\mathcal{C}(G)$, is a vector space over ${\mathbb F}_2$, spanned by all incidence vectors of edge-sets of cycles of $G$. If $G$ has $n$ vertices, then $\mathcal{C}_n(G)$ is the subspace of $\mathcal{C}(G)$, spanned by the incidence vectors of Hamilton cycles of $G$. We prove that asymptotically almost surely $\mathcal{C}_n(\mathcal{G}_{n,d}) = \mathcal{C}(\mathcal{G}_{n,d})$ holds whenever $n$ is odd and $d$ is a sufficiently large (even) integer. This extends (though with a weaker bound on $d$) the well-known result asserting that $\mathcal{G}_{n,d}$ is asymptotically almost surely Hamiltonian for every $d \geq 3$ (but not for $d < 3$). Since $n$ being odd mandates that $d$ be even, somewhat limiting the generality of our result, we also prove that if $n$ is even and $d$ is any sufficiently large integer, then asymptotically almost surely $\mathcal{C}_{n-1}(\mathcal{G}_{n,d}) = \mathcal{C}(\mathcal{G}_{n,d})$.

An influential result of Bohman, Frieze, and Martin asserts that if $H$ is an $n$-vertex graph with minimum degree at least $\delta n$ for some constant $\delta > 0$, and $G \sim \mathbb{G}(n, C/n)$, where $C := C(\delta)$ is a sufficiently large constant, then $H \cup G$ is asymptotically almost surely Hamiltonian. We strengthen this result by proving that the same assumptions on $H$ and $G$ ensure that $\mathcal{C}_n(H \cup G) = \mathcal{C}(H \cup G)$ holds asymptotically almost surely.        
\end{abstract}

\section{Introduction}  
Let $G = (V,E)$ be a graph on $n$ vertices. The \emph{edge space} of $G$, denoted $\mathcal{E}(G)$, is a vector space over ${\mathbb F}_2$ consisting of all incidence vectors of subsets of $E$. The \emph{cycle space} of $G$, denoted $\mathcal{C}(G)$, is the subspace of $\mathcal{E}(G)$, spanned by all incidence vectors of cycles of $G$. For any integer $3 \leq k \leq n$, let $\mathcal{C}_k(G)$ be the subspace of $\mathcal{C}(G)$, spanned by all incidence vectors of cycles of length $k$ in $G$. Determining conditions under which $\mathcal{C}_k(G) = \mathcal{C}(G)$ holds for some $3 \leq k \leq n$ is a well-studied problem (see, e.g.,~\cite{BK, BL, DHJ, Hartman, Locke1, Locke2}). In this paper we are interested in the case $k = n$, that is, in graphs whose cycle space is spanned by their Hamilton cycles. This problem has been addressed by various researchers (see, e.g.,~\cite{ALW, HK, Heinig1, Heinig2, HY}). Since the symmetric difference of any two even graphs (i.e., subsets of $E$ of even size) is an even graph, it is evident that if $\mathcal{C}_n(G) = \mathcal{C}(G)$, then $G$ is bipartite or $n$ is odd. Moreover, $G$ must either be acyclic or Hamiltonian. Since the former case is not very interesting, this study can be viewed as part of the common theme of proving that (possibly slightly strengthened) various sufficient conditions for Hamiltonicity in fact ensure stronger properties. 

A natural venue for this problem are random and pseudo-random graphs. Indeed, it was proved by Christoph, Nenadov, and Petrova in~\cite{CNP} that $\mathcal{C}_n(G) = \mathcal{C}(G)$ holds whenever $G$ is a pseudo-random graph with certain appropriate properties (see Theorem 1.3 in~\cite{CNP} for details). One immediate consequence of this result is that if $G \sim \mathbb{G}(n,p)$, where $n$ is odd and $p \geq C \log n/n$ for a sufficiently large constant $C$, then asymptotically almost surely (a.a.s. hereafter) $\mathcal{C}_n(G) = \mathcal{C}(G)$. The exact bound on $p$ has been subsequently attained by the authors in~\cite{HK} where it is proved that if $p$ is large enough to ensure that $\delta(G) \geq 3$ holds a.a.s., then it also ensures that a.a.s. $\mathcal{C}_n(G) = \mathcal{C}(G)$ (the necessity of this condition was observed by Heinig in~\cite{Heinig2}).   

Another immediate consequence of the aforementioned result of Christoph, Nenadov, and Petrova regarding pseudo-random graphs is that if $G$ is an $(n, d, \lambda)$-graph, $n$ is odd, $d \geq C \log n$, and $\lambda \leq \varepsilon d/\log n$, where $C$ and $\varepsilon$ are appropriate constants, then $\mathcal{C}_n(G) = \mathcal{C}(G)$. We briefly discuss $(n, d, \lambda)$-graphs in Section~\ref{sec::openprob} where we improve this result. However, in this paper our main focus is on random regular graphs $G \sim \mathcal{G}_{n, d}$ for which we prove that sufficiently large yet constant $d$ is a.a.s. enough to ensure that $\mathcal{C}_n(G) = \mathcal{C}(G)$ holds.
 
\begin{theorem} \label{th::CycleSpaceGnd}
There exists an integer $d_0$ such that the following holds for every $d \geq d_0$. Let $G \sim \mathcal{G}_{n, d}$, where $n$ is odd. Then, a.a.s. $\mathcal{C}_n(G) = \mathcal{C}(G)$. 
\end{theorem}

While, formally, the assertion of Theorem~\ref{th::CycleSpaceGnd} holds for every sufficiently large $d$, the required assumption that $n$ is odd  immediately implies that $d$ is even. The following result constitutes a remedy to this limitation.

\begin{theorem} \label{th::evenn}
There exists an integer $d_0$ such that the following holds for every $d \geq d_0$. Let $G \sim \mathcal{G}_{n, d}$, where $n$ is even. Then, a.a.s. $\mathcal{C}_{n-1}(G) = \mathcal{C}(G)$. 
\end{theorem}

Another natural venue for the problem of determining whether $\mathcal{C}_n(G) = \mathcal{C}(G)$ holds for an $n$-vertex Hamiltonian graph $G$ is that of Dirac graphs. A classical result of Dirac~\cite{Dirac} asserts that every $n$-vertex graph $H$ (on at least three vertices) with minimum degree $\delta(H) \geq n/2$ is Hamiltonian. As for its cycle space, it was proved by Heinig~\cite{Heinig1} that $\mathcal{C}_n(G) = \mathcal{C}(G)$ holds whenever $G$ is an $n$-vertex graph with minimum degree $\delta(G) \geq (1/2 + \varepsilon) n$, for any $\varepsilon > 0$ and sufficiently large odd $n$. This was significantly improved by Christoph, Nenadov, and Petrova~\cite{CNP} who proved that $\delta(G) \geq n/2 + 41$ suffices. Recently, this problem has been completely resolved by Hou and Yin~\cite{HY} who proved that $\delta(G) > n/2$ is sufficient (and necessary) for $\mathcal{C}_n(G) = \mathcal{C}(G)$ to hold. In light of this result and since, as noted above, the corresponding problem for $\mathbb{G}(n,p)$ has been resolved in~\cite{HK}, it seems natural to consider this cycle space problem for randomly perturbed dense graphs.

The study of randomly perturbed graphs was initiated by Bohman, Frieze, and Martin~\cite{BFM} who discovered that once slightly randomly perturbed, graphs with linear minimum degree become Hamiltonian asymptotically almost surely. Formally, they proved that for every constant $\delta > 0$ there exists a constant $C := C(\delta)$ such that $H \cup R$ is a.a.s. Hamiltonian, whenever $H$ is an $n$-vertex graph with minimum degree at least $\delta n$ and $R \sim \mathbb{G}(n, C/n)$, undershooting the threshold for Hamiltonicity in $\mathbb{G}(n, p)$ by a logarithmic factor. Such a result can
be seen as bridging between results regarding the Hamiltonicity of dense graphs and the emergence of such cycles in random graphs. Various extensions of this result have subsequently been proved (see, e.g.,~\cite{AHK, HMMMP}). Our second main result asserts that a.a.s. the Hamilton cycles of a randomly perturbed graph, as in the aforementioned result of Bohman, Frieze, and Martin, span its cycle space.
\begin{theorem} \label{th::CycleSpacePerturbed}
For every constant $\delta > 0$ there exists a constant $C := C(\delta)$ such that the following holds. Let $H$ be an $n$-vertex graph, where $n$ is odd, with minimum degree $\delta(H) \geq \delta n$, and let $G \sim \mathbb{G}(n,p)$, where $p := p(n) \geq C/n$. Then a.a.s. $\mathcal{C}_n(H \cup G) = \mathcal{C}(H \cup G)$. 
\end{theorem}

The rest of this paper is organised as follows. In Section~\ref{sec::prelim} we introduce some terminology, notation, and standard tools, and present the method of Christoph, Nenadov, and Petrova from~\cite{CNP} which is a central ingredient in our proofs. Starting with the result with the simpler proof, we prove Theorem~\ref{th::CycleSpacePerturbed} in Section~\ref{sec::perturbed}. In Section~\ref{sec::random} we prove Theorem~\ref{th::CycleSpaceGnd}, and then briefly explain how to adjust our proof so as to obtain a proof of Theorem~\ref{th::evenn}. Finally, in Section~\ref{sec::openprob} we present several open problems.

\section{Preliminaries and tools} \label{sec::prelim}

For the sake of simplicity and clarity of presentation, we do not make a particular effort to optimize the constants obtained in some of our proofs. We also omit floor and ceiling signs whenever these are not crucial. Most of our results are asymptotic in nature and whenever necessary we assume that the number of vertices $n$ is sufficiently large. Throughout this paper, $\log$ stands for the natural logarithm, unless explicitly stated otherwise. Our graph-theoretic notation is standard; in particular, we use the following.

For a graph $G$, let $V(G)$ and $E(G)$ denote its sets of vertices and edges respectively, and let $v(G) = |V(G)|$ and $e(G) = |E(G)|$. For a set $A \subseteq V(G)$, let $E_G(A)$ denote the set of edges of $G$ with both endpoints in $A$ and let $e_G(A) = |E_G(A)|$. For disjoint sets $A, B \subseteq  V(G)$, let $E_G(A, B)$ denote the set of edges of $G$ with one endpoint in $A$ and one endpoint in $B$, and let $e_G(A, B) = |E_G(A, B)|$. For a set $S \subseteq V(G)$, let $G[S]$ denote the subgraph of $G$ induced by the set $S$. For a set $S \subseteq V(G)$, let $N_G(S) = \{v \in V(G) \setminus S : \exists u \in S \textrm{ such that } uv \in E(G)\}$ denote the \emph{external neighbourhood} of $S$ in $G$. For a vertex $u \in V(G)$ we abbreviate $N_G(\{u\})$ under $N_G(u)$ and let $\textrm{deg}_G(u) = |N_G(u)|$ denote the degree of $u$ in $G$. The maximum degree of a graph $G$ is $\Delta(G) = \max\{\textrm{deg}_G(u) : u \in V(G)\}$, and the minimum degree of a graph $G$ is $\delta(G) = \min \{\textrm{deg}_G(u) : u \in V(G)\}$. For a vertex $u \in V(G)$ and a set $S \subseteq V(G)$, let $\textrm{deg}_G(u, S) = |N_G(u) \cap S|$. For a vertex $x \in V(G)$, let $\partial_G(x) = \{xy : y \in N_G(x)\}$. Given any two (not necessarily distinct) vertices $x, y \in V(G)$, the \emph{distance} between $x$ and $y$ in $G$, denoted $\textrm{dist}_G(x,y)$, is the length of a shortest path between $x$ and $y$ in $G$, where the length of a path is the number of its edges (for the sake of formality, we define $\textrm{dist}_G(x,y)$ to be $\infty$ whenever $x$ and $y$ lie in different connected components of $G$). The \emph{diameter} of $G$, denoted $\textrm{diam}(G)$, is $\max \{\textrm{dist}_G(x,y) : x, y \in V(G)\}$. 

Let $n$ be an odd integer, and let $G$ be an $n$-vertex Hamiltonian graph. A recipe for proving $\mathcal{C}_n(G) = \mathcal{C}(G)$ is presented in~\cite{CNP}. In order to describe it we need some definitions and results. 

\begin{lemma} [\cite{CNP}] \label{lem::subgraphR}
Let $G$ be an $n$-vertex Hamiltonian graph, where $n$ is odd, and suppose that $\mathcal{C}_n(G) \neq \mathcal{C}(G)$. Then, there
exists a subgraph $R$ of $G$ such that the following conditions hold.
\begin{enumerate}
\item [\emph{(C1)}] $R \neq G$;

\item [\emph{(C2)}] Every Hamilton cycle in $G$ contains an even number of edges from $R$;

\item [\emph{(C3)}] For every partition $V (G) = A \cup B$ it holds that $e_R(A, B) \geq e_G(A, B)/2$ and $R \neq G[A, B]$.
\end{enumerate}
\end{lemma}

The following definition of a so-called \emph{parity switcher} is central to the method of~\cite{CNP}. It describes a construction that, given graphs $G$ and $R$ as in Lemma~\ref{lem::subgraphR}, aids one in finding a Hamilton cycle of $G$ with an odd number of edges in $R$, thus arriving at a contradiction to (C2) above. 

\begin{definition} \label{def::paritySwitcher}
Given a graph $G$ and a subgraph $R \subseteq G$, a subgraph $W \subseteq G$ is called an $R$-\emph{parity-switcher} if it consists of an even cycle $C = (v_1, v_2, \ldots, v_{2k}, v_1)$ with an odd number of edges in $R$, and vertex-disjoint paths $P_2, \ldots, P_k$ such that $\bigcup_{i=2}^k E(P_i) \cap E(C) = \varnothing$ and, for every $2 \leq i \leq k$, the endpoints of $P_i$ are $v_i$ and $v_{2k-i+2}$.
\end{definition}

We may now specify the recipe from~\cite{CNP}; it consists of the following five steps.
\begin{enumerate}

\item [(S1)] Let $G$ be an $n$-vertex Hamiltonian graph, where $n$ is odd. Suppose it satisfies $\mathcal{C}_n(G) \neq \mathcal{C}(G)$, and let $R \subseteq G$ be a subgraph as in Lemma~\ref{lem::subgraphR}.

\item [(S2)] Find in $G$ a (small) $R$-parity-switcher $W$, that is,
\begin{enumerate}
\item [(S2a)] Find an even (short) cycle $C = (v_1, \ldots, v_{2k}, v_1)$ with an odd number of edges in $R$.

\item [(S2b)] Find pairwise vertex-disjoint (short) paths $P_i$ between $v_i$ and $v_{2k-i+2}$ for every $2 \leq i \leq k$.
\end{enumerate}

\item [(S3)] Find in $(G \setminus V(W)) \cup \{v_1, v_{k+1}\}$ a Hamilton path $P$ whose endpoints are $v_1$ and $v_{k+1}$.

\item [(S4)] If $P$ contains an odd (even) number of edges of $R$, then choose a Hamilton path $P'$ in $W$ whose endpoints are $v_1$ and $v_{k+1}$  with an even (odd) number of edges of $R$.

\item [(S5)] Conclude that the concatenation of $P$ and $P'$ yields a Hamilton cycle $H \subseteq G$ with an odd number of edges in R, contradicting (C2).

\end{enumerate}

Note that there is nothing to prove in steps (S4) and (S5). Moreover, whenever we start with a graph which we know to be Hamiltonian, step (S1) becomes immediate. The main task is thus to deal with steps (S2) and (S3). 

\bigskip

The following result from~\cite{CNP} is our main tool for handling Step (S2a).
\begin{lemma} \label{lem::S2CNP}
Let $R \subseteq G$ be graphs, and let $\ell$ be a positive integer. If
\begin{enumerate}
\item [\emph{(L1)}] For every $S \subseteq V(G)$ of size $|S| \leq 2 \ell$ and any two vertices $x, y \in V(G) \setminus S$, there is a path between $x$
and $y$ in $R \setminus S$ whose length is at most $\ell - 1$, and

\item [\emph{(L2)}] $R \neq G$ and $R \neq G[A, B]$ for every partition $V(G) = A \cup B$,
\end{enumerate}
then there exists an even cycle $C \subseteq G$ of length $|C| \leq 2\ell$ that contains an odd number of edges from $R$.
\end{lemma}

The following known result (see Theorem~\ref{th::HamConCexpander} below) is our main tool for handling Step (S3). In order to state it, we require the notion of an \emph{expander} and the notion of \emph{Hamilton-connectivity}.

\begin{definition} \label{def::Cexpander}
An $n$-vertex graph $G$, where $n \geq 3$, is called a $c$-expander if it satisfies the following two properties.
\begin{enumerate}
\item [\emph{(E1)}] $|N_G(X)| \geq c |X|$ holds for every $X \subseteq V(G)$ of size $|X| < n/(2c)$;

\item [\emph{(E2)}] There is an edge of $G$ between any two disjoint sets $X, Y \subseteq V(G)$ of size $|X|, |Y | \geq n/(2c)$.
\end{enumerate}
\end{definition}

\begin{definition} \label{def::HanCon}
A graph $G$ is said to be Hamilton-connected if for every two vertices $x, y \in V(G)$ there is a Hamilton path of $G$ whose endpoints are $x$ and $y$.
\end{definition}

\begin{theorem} [Theorem 7.1 in~\cite{DMMPS}] \label{th::HamConCexpander}
For every sufficiently large $c > 0$, every $c$-expander is Hamilton-connected.
\end{theorem}

\section{Randomly perturbed dense graphs} \label{sec::perturbed}

The main aim of this section is to prove Theorem~\ref{th::CycleSpacePerturbed}. Before doing so, we state and prove several auxiliary results that will facilitate our proof. 

\bigskip

The following known result is our main tool for handling Step (S1) (it is a rephrasing of part of a result from~\cite{BFM}).
\begin{theorem} [Theorem 1 in~\cite{BFM}] \label{th::perturbedHam}
For every constant $\delta > 0$ there exists a constant $C := C(\delta)$ such that the following holds. Let $H$ be an $n$-vertex graph with minimum degree $\delta(H) \geq \delta n$ and let $G \sim \mathbb{G}(n,p)$, where $p \geq C/n$. Then, $H \cup G$ is a.a.s. Hamiltonian.
\end{theorem}

Next, we state and prove several properties of dense graphs (acting as the seed, i.e., the graph being perturbed) and of sparse random graphs (acting as the random perturbation).
\begin{lemma} \label{lem::randomGraph}
For all constants $0 < \alpha, \beta < 1$ and $c \geq 0$ there exists a constant $C := C(\alpha, \beta, c)$ such that the following holds. Let $G \sim \mathbb{G}(n,p)$, where $p \geq C/n$. Then, a.a.s. $e_G(A,B) > c n$ holds for every $A \subseteq [n]$ of size $|A| = \alpha n$ and every $B \subseteq [n] \setminus A$ of size $|B| = \beta n$.
\end{lemma}

\begin{proof}
Fix any two sets $A \subseteq [n]$ of size $|A| = \alpha n$ and $B \subseteq [n] \setminus A$ of size $|B| = \beta n$. Note that $e_G(A,B) \sim \textrm{Bin}(|A| |B|, p)$ and thus, in particular, $\mathbb{E}(e_G(A,B)) = |A| |B| p$. Then, 
$$
\mathbb{P}(e_G(A,B) \leq c n) \leq \mathbb{P}(e_G(A,B) \leq \mathbb{E}(e_G(A,B))/2) \leq \exp \left\{- \frac{C}{8 n} \cdot \alpha n \cdot \beta n \right\} = o(4^{- n}),
$$
where the first inequality and the equality hold by an appropriate choice of $C$, and the second inequality holds by Chernoff's bound. The assertion of the lemma then follows by a union bound over all choices of $A$ and $B$.
\end{proof}

\begin{lemma} \label{lem::robustExpansion}
For every constant $\delta > 0$ and every sufficiently large constant $c$, there exists a constant $C := C(\delta, c)$ such that the following holds. Let $H$ be an $n$-vertex graph with minimum degree $\delta(H) \geq \delta n$, and let $G \sim \mathbb{G}(n,p)$, where $p := p(n) \geq C/n$. Then, a.a.s. $(H \cup G)[[n] \setminus U]$ is a $c$-expander for every $U \subseteq [n]$ of size $|U| \leq \delta n/10$.
\end{lemma}

\begin{proof}
Let $\alpha = \delta/(3c)$, let $\beta = 1/(3c)$, and let $C := C(\alpha, \beta, 0)$ be as in the statement of Lemma~\ref{lem::randomGraph}. It then follows by Lemma~\ref{lem::randomGraph} that a.a.s. $E_G(A,B) \neq \varnothing$ holds for every $A \subseteq [n]$ of size $|A| = \delta n/(3c)$ and every $B \subseteq [n] \setminus A$ of size $|B| = n/(3c)$. In particular, this implies Property (E2) from Definition~\ref{def::Cexpander} for $(H \cup G)[[n] \setminus U]$ and any set $U \subseteq [n]$ of size $|U| \leq \delta n/10$.

Next, we prove that a.a.s. $(H \cup G)[[n] \setminus U]$ satisfies Property (E1) from Definition~\ref{def::Cexpander} for every $U \subseteq [n]$ of size $|U| \leq \delta n/10$. Fix an arbitrary set $U \subseteq [n]$ of size $|U| \leq \delta n/10$. Fix an arbitrary set $X \subseteq [n] \setminus U$ of size $|X| \leq \frac{n}{2c}$. Assume first that $|X| \leq \frac{\delta n}{2c}$; we may assume that $X \neq \varnothing$. Then,
$$
|N_{H \cup G}(X) \cap ([n] \setminus U)| \geq |N_H(X) \setminus U| \geq \delta(H) - |X| - \delta n/10 \geq \delta n - \delta n/(2c) - \delta n/10 \geq \delta n/2 \geq c |X|.
$$ 
Assume then that $\frac{\delta n}{2c} \leq |X| \leq \frac{n}{2c}$. If $|N_{H \cup G}(X) \cap ([n] \setminus U)| < c |X|$, then there exist disjoint sets $A \subseteq X$ and $B \subseteq [n] \setminus (X \cup (N_{H \cup G}(X) \cap ([n] \setminus U)) \cup U)$ such that $|A| = \delta n/(3c)$, $|B| = n/(3c)$, and $E_G(A,B) = \varnothing$. However, as noted in the previous paragraph, the probability that there exist such sets $A$ and $B$ is $o(1)$.
\end{proof}


\begin{lemma} \label{lem::smallDiameter}
Let $H = (V,E)$ be a connected $n$-vertex graph with minimum degree $d$. Then, $\emph{diam}(H) < 3n/d$.
\end{lemma}

\begin{proof}
Fix any two vertices $u, v \in V$. Since $H$ is connected, there exists a path between $u$ and $v$ in $H$; let $u = x_1, x_2, \ldots, x_t = v$ be a shortest such path, and suppose for a contradiction that $t \geq 3n/d + 1$. Then, for every $0 \leq i < j \leq n/d$, it must hold that $N_H(x_{3i+1}) \cap N_H(x_{3j+1}) = \varnothing$. It thus follows that 
$$
\left|\bigcup_{i=0}^{n/d} N_H(x_{3i+1}) \right| = \sum_{i=0}^{n/d} |N_H(x_{3i+1})| \geq d (n/d + 1) > n,
$$
which is an obvious contradiction.
\end{proof}

\begin{lemma} \label{lem::robustConnectivity}
Let $H$ be an $n$-vertex graph, where $n$ is sufficiently large, with minimum degree $\delta(H) \geq \delta n$ for some constant  $\delta > 0$, and let $S \subseteq V(H)$ be a set of size $r = o(n)$. Suppose that $e_H(A,B) > r n$ holds for every $A \subseteq V(H)$ of size $|A| \geq \delta n/3$ and every $B \subseteq V(H) \setminus A$ of size $|B| \geq n/3$. Let $R$ be a subgraph of $H$ such that $e_R(A,B) \geq e_H(A,B)/2$ for every partition $V(H) = A \cup B$. Then, $R \setminus S$ is connected.
\end{lemma}

\begin{proof}
Suppose for a contradiction that there exists a set $S \subseteq V(H)$ of size $r$ such that $R \setminus S$ is disconnected. Let $X$ be the vertex-set of a connected component of $R \setminus S$ and let $Y = V(H) \setminus (S \cup X)$; assume without loss of generality that $|X| \leq |Y|$. Since $\delta(H) \geq \delta n$ and $e_R(A,B) \geq e_H(A,B)/2$ for every partition $V(H) = A \cup B$, it follows that $\delta(R) \geq \delta n/2$. Since, moreover, $r = o(n)$, it must hold that $|X| \geq \delta n/3$ and $|Y| \geq n/3$. Therefore
$$
0 = e_{R \setminus S}(X, Y) \geq e_R(X, Y \cup S) - |S| |X| \geq e_H(X, Y \cup S)/2 - r n/2 > r n/2 - r n/2 = 0,
$$
which is an obvious contradiction.
\end{proof}

\begin{proof} [Proof of Theorem~\ref{th::CycleSpacePerturbed}]
Let $c$ be a sufficiently large constant as per Theorem~\ref{th::HamConCexpander}. Let $\alpha = \delta/(3c)$, let $\beta = 1/(3c)$, and let $r = 100 \delta^{-2}$. Let $C_1 := C_1(\alpha, \beta, r)$ be the constant whose existence is ensured by Lemma~\ref{lem::randomGraph}, let $C_2 := C_2(\delta, c)$ be the constant whose existence is ensured by Lemma~\ref{lem::robustExpansion}, let $C_3 := C_3(\delta)$ be the constant whose existence is ensured by Theorem~\ref{th::perturbedHam}, and let $C = \max \{C_1, C_2, C_3\}$. Let $H$ and $G$ (with $C$ as above) be as in the premise of the theorem; by our choice of $C$ we may then assume that $H$ and $G$ satisfy the assertions of Theorem~\ref{th::perturbedHam} and of Lemmas~\ref{lem::randomGraph} and~\ref{lem::robustExpansion}. Suppose for a contradiction that $\mathcal{C}_n(H \cup G) \neq \mathcal{C}(H \cup G)$. Let $R$ be a subgraph of $H \cup G$ as per Lemma~\ref{lem::subgraphR}. It then follows by Theorem~\ref{th::perturbedHam} and by Lemmas~\ref{lem::subgraphR}, \ref{lem::randomGraph}, \ref{lem::robustExpansion}, \ref{lem::smallDiameter}, and \ref{lem::robustConnectivity} that $H \cup G$ and $R$ satisfy all of the following properties:

\begin{enumerate}

\item [(i)] $H \cup G$ is Hamiltonian;

\item [(ii)] $(H \cup G)[[n] \setminus U]$ is a $c$-expander for every $U \subseteq [n]$ of size $|U| \leq \delta n/10$;

\item [(iii)] $\delta(R) \geq \delta(H \cup G)/2 \geq \delta(H)/2 \geq \delta n/2$;

\item [(iv)] $\textrm{diam}(R \setminus U) < 7/\delta$ for every $U \subseteq [n]$ of size $|U| \leq r$.

\end{enumerate}

We follow the recipe that was presented in Section~\ref{sec::prelim}. That is, we need to handle steps (S1), (S2), and (S3). Our assumption that $\mathcal{C}_n(H \cup G) \neq \mathcal{C}(H \cup G)$ will then lead to the contradiction appearing in (S5).  

Starting with (S1), we note that it is an immediate corollary of Property (i) above, and our choice of $R$.

Next, we take care of (S2). Starting with (S2a), by properties (C1) and (C3) of Lemma~\ref{lem::subgraphR} the graph $R$ satisfies Property (L2) of Lemma~\ref{lem::S2CNP}. Moreover, by Property (iv) above, $R$ satisfies Property (L1) of Lemma~\ref{lem::S2CNP} with $\ell := 7/\delta$. It thus follows by Lemma~\ref{lem::S2CNP} that $G$ contains an even cycle $C = (v_1, \ldots, v_{2k})$, for some $k \leq \ell$, having an odd number of edges in $R$. As for (S2b), we construct the required paths $P_2, \ldots, P_k$ one by one as follows. Assume that for some $2 \leq i \leq k$ we have already built $P_2, \ldots, P_{i-1}$, each of length at most $\ell$, and now wish to construct $P_i$. Let $W_i = \left(\{v_1, \ldots, v_{2k}\} \cup V(P_2) \cup \ldots \cup V(P_{i-1}) \right) \setminus \{v_i, v_{2k-i+2}\}$ and note that $|W_i| \leq \ell^2 \leq r$. By Property (iv) above, we conclude that there is a path $P_i$ of length at most $\ell$ in $(H \cup G) \setminus W_i$ between $v_i$ and $v_{2k-i+2}$.
 
Finally, we establish (S3). Let $W = \left(\{v_1, \ldots, v_{2k}\} \cup V(P_2) \cup \ldots \cup V(P_k) \right) \setminus \{v_1, v_{k+1}\}$, and note that $|W| \leq 2 k + (k-1) (\ell - 1) - 2 \leq \ell^2 \leq r = o(n)$. It follows by Property (ii) above that $(H \cup G)[[n] \setminus W]$ is a $c$-expander. It then follows by Theorem~\ref{th::HamConCexpander} that $(H \cup G)[[n] \setminus W]$ admits a Hamilton path whose endpoints are $v_1$ and $v_{k+1}$.
\end{proof}

\section{Random regular graphs} \label{sec::random}

The main aim of this section is to prove Theorem~\ref{th::CycleSpaceGnd}. Before doing so, we state and prove several auxiliary results that will facilitate our proof; some are well-known and some are new. 

As already noted in the introduction, $\mathcal{G}_{n, d}$ is a.a.s. Hamiltonian for almost any value of $d$.
\begin{theorem} [\cite{RW1, RW2}] \label{th::HamGnd}
Let $G \sim \mathcal{G}_{n,d}$, where $d \geq 3$. Then $G$ is a.a.s. Hamiltonian.
\end{theorem}

The following lemma is an immediate consequence of the expander mixing lemma (see, e.g.,~\cite{AS}) and the fact that $\mathcal{G}_{n,d}$ is an $(n, d, \lambda)$-graph such that a.a.s. $\lambda(\mathcal{G}_{n,d}) \leq 2 \sqrt{d - 1} + \varepsilon$, where $\varepsilon > 0$ is arbitrarily small yet fixed~\cite{Friedman}.

\begin{lemma} \label{lem::edgeDistribution}
Let $G \sim \mathcal{G}_{n,d}$, where $d$ is a sufficiently large constant. Then a.a.s. the following properties hold
\begin{enumerate}
\item [$(a)$] $\left|e_G(A,B) - \frac{|A| |B| d}{n} \right| \leq 3 \sqrt{d |A| |B|}$ holds for any two disjoint sets $A, B \subseteq V(G)$;

\item [$(b)$] $\left|e_G(A) - \frac{d}{n} \binom{|A|}{2} \right| \leq 3 \sqrt{d} |A|$ holds for every set $A \subseteq V(G)$.
\end{enumerate}
\end{lemma}

\begin{lemma} \label{lem::halfDegree}
Let $\varepsilon > 0$ be an arbitrarily small constant, and let $G \sim \mathcal{G}_{n,d}$, where $d$ is a sufficiently large constant. Then, a.a.s. for every $A \subseteq V(G)$ of size $|A| \geq \varepsilon n$ and every $B \subseteq V(G) \setminus A$ of size $|B| = 0.59 n$ there exists a vertex $u \in A$ such that $\emph{deg}_G(u, B) \geq 0.51 d$.
\end{lemma}

\begin{proof}
Observe that it suffices to prove the lemma for all sets $A$ whose size is precisely $\varepsilon n$. It follows by Lemma~\ref{lem::edgeDistribution}(a) that a.a.s. for every $A \subseteq V(G)$ of size $\varepsilon n$ and every $B \subseteq V(G) \setminus A$ of size $|B| = 0.59 n$ it holds that
$$
e_G(A, B) \geq 0.58 \varepsilon n d.
$$
Assume then that $G$ satisfies this property, and suppose for a contradiction that there exist sets $A \subseteq V(G)$ of size $\varepsilon n$ and $B \subseteq V(G) \setminus A$ of size $|B| = 0.59 n$ such that $\textrm{deg}_G(u, B) < 0.51 d$ holds for every $u \in A$. It then follows that
\begin{align*}
 0.58 \varepsilon n d \leq e_G(A, B) = \sum_{u \in A} \textrm{deg}_G(u, B) < \sum_{u \in A} 0.51 d = 0.51 \varepsilon n d,
\end{align*}
which is a clear contradiction.
\end{proof}

The following result allows us to split a graph into several parts in a beneficial manner. The specific formulation we use is taken from~\cite{HKT}, though similar results can be found in other sources.
\begin{lemma} [Lemma 2.4 in~\cite{HKT}] \label{lem::LLLsplit}
Let $G = (V, E)$ be a graph on $n$ vertices with maximum degree $\Delta$. Let $Y \subseteq V$ be a set of $m = a + b$ vertices, where $a$ and $b$ are positive integers. Assume that $\emph{deg}_G(v, Y) \geq \delta$ holds for every $v \in V$. If $\Delta^2 \cdot \lceil \frac{m}{\min \{a,b\}} \rceil \cdot 2 \cdot e^{1 - \frac{\min \{a,b\}^2}{5 m^2} \delta} < 1$, then there exists a partition $Y = A \cup B$ of $Y$ such that the following properties hold.
\begin{enumerate}
\item [$(1)$] $|A| = a$ and $|B| = b$;

\item [$(2)$] $\emph{deg}_G(v, A) \geq \frac{a}{3m} \emph{deg}_G(v, Y)$ holds for every $v \in V$;

\item [$(3)$] $\emph{deg}_G(v, B) \geq \frac{b}{3m} \emph{deg}_G(v, Y)$ holds for every $v \in V$.  
\end{enumerate}
\end{lemma}

\begin{lemma} \label{lem::subgraphExpander}
For every $\alpha \in [0,1)$ and $\delta \in (0,1]$ there exists $c := c(\alpha, \delta)$ such that the following holds. Let $G \sim \mathcal{G}_{n,d}$, where $d$ is a sufficiently large constant. Let $H$ be a (not necessarily spanning) subgraph of $G$ satisfying $\delta(H) \geq \delta d$. Then, a.a.s. $|N_{H \setminus F}(X)| \geq c |X| d$ for every $F \subseteq E(H)$ and every $X \subseteq V(H)$ of size $|X| \leq n/d$ such that $\left|F \cap \partial_H(x) \right| \leq \alpha \cdot \emph{deg}_H(x)$ holds for every $x \in X$. 
\end{lemma}

\begin{proof}
Suppose that $X \subseteq V(H)$ is a set of size $|X| \leq n/d$ and $F \subseteq E(H)$ is such that $\left|F \cap \partial_H(x) \right| \leq \alpha \cdot \textrm{deg}_H(x)$ holds for every $x \in X$, and yet $|N_{H \setminus F}(X)| < c |X| d$. It follows by the premise of the lemma and by Lemma~\ref{lem::edgeDistribution} that a.a.s.
\begin{align*}
|X| \cdot (1 - \alpha) \delta d &\leq \sum_{x \in X} \textrm{deg}_{H \setminus F}(x) \leq 2 e_G(X) + e_G(X, N_{H \setminus F}(X)) \\
&\leq \left(\frac{d |X|^2}{n} + 6 \sqrt{d} |X| \right) + \left(\frac{d |X| \cdot c |X| d}{n} + 3 \sqrt{c} d |X| \right) \\
&\leq |X| (1 + 6 \sqrt{d} + c d + 3 \sqrt{c} d), 
\end{align*}
which is a contradiction by our choice of $d$ being sufficiently large and $c$ being sufficiently small with respect to $\delta$ and $1 - \alpha$. 
\end{proof}

\begin{lemma} \label{lem::robustEdgeBetweenSets}
Let $G \sim \mathcal{G}_{n,d}$, where $d$ is a sufficiently large constant. Then a.a.s. the following holds. Let $R$ be a subgraph of $G$ such that $e_R(A, V(G) \setminus A) \geq e_G(A, V(G) \setminus A)/2$ for every $A \subseteq V(G)$. Then, $e_R(A, B) > 0$ holds for any two disjoint sets $A, B \subseteq V(G)$ of size $|A| = |B| = 2n/5$. 
\end{lemma}

\begin{proof}
Fix any two disjoint sets $A, B \subseteq V(G)$ of size $|A| = |B| = 2n/5$, and let $D = V(G) \setminus (A \cup B)$. It follows by Lemma~\ref{lem::edgeDistribution}(a) that a.a.s. 
\begin{align} \label{eq::eAnotA}
e_G(A, V(G) \setminus A) &\geq d |A| |V(G) \setminus A|/n - 3 \sqrt{d |A| |V(G) \setminus A|} \geq \frac{6}{25} \cdot d n - 2 \sqrt{d} n,
\end{align} 
and that 
\begin{align} \label{eq::eAD}
e_G(A, D) &\leq d |A| |D|/n + 3 \sqrt{d |A| |D|} \leq \frac{2}{25} \cdot d n + \sqrt{d} n. 
\end{align} 
Combining~\eqref{eq::eAnotA} and~\eqref{eq::eAD} we obtain
\begin{align*}
e_R(A, B) &= e_R(A, V(G) \setminus A) - e_R(A,D) \geq e_G(A, V(G) \setminus A)/2 - e_G(A,D) \\
&\geq \left(\frac{3}{25} \cdot d n - \sqrt{d} n \right) - \left(\frac{2}{25} \cdot d n + \sqrt{d} n \right) > 0,
\end{align*}
where the last inequality holds for sufficiently large $d$.
\end{proof}

\begin{lemma} \label{lem::DiameterR}
Let $G \sim \mathcal{G}_{n,d}$, where $d$ is a sufficiently large constant. Then a.a.s. the following holds. Let $R$ be a subgraph of $G$ with minimum degree $\delta(R) \geq d/2$. Then, there exists an absolute constant $K$ such that for every set $S \subseteq V(G)$ of size $|S| = o(n)$ and such that $\emph{deg}_G(u, S) \leq d/100$ holds for every $u \in V(G)$, and for every two vertices $x, y \in V(G) \setminus S$, there is a path between $x$ and $y$ in $R \setminus S$ whose length is at most $\frac{K \log n}{\log d}$.
\end{lemma}

\begin{proof}
Fix an arbitrary set $S \subseteq V(G)$ of size $|S| = o(n)$ and such that $\textrm{deg}_G(u, S) \leq d/100$ holds for every $u \in V(G)$. Since, moreover, $\delta(R) \geq d/2$ holds by the premise of the lemma, it follows that $\textrm{deg}_{R \setminus S}(x) \geq 0.49 d$ holds for every $x \in V(G) \setminus S$. For every vertex $x \in V(G) \setminus S$ and every non-negative integer $i$, let $N_{R \setminus S}^i(x) = \{u \in V(G) \setminus S : \textrm{dist}_{R \setminus S}(x,u) \leq i\}$. Starting from an arbitrary vertex $x \in V(G) \setminus S$, repeated applications of Lemma~\ref{lem::subgraphExpander} with $\alpha = 0$ and $H = R \setminus S$ show that a.a.s. there exists an integer $t \leq \frac{\log n}{\log (c d)}$, where $c := c(0, 0.49)$ is the constant whose existence is ensured by Lemma~\ref{lem::subgraphExpander}, such that $\left|N_{R \setminus S}^t(x) \right| \geq c n$. We claim that a.a.s. $\left|N_{R \setminus S}^{t+1}(x) \right| \geq 2 n/5$. Indeed, if not, then there exists a set $B \subseteq V(G) \setminus \left(S \cup N_{R \setminus S}^{t+1}(x) \right)$ of size $|B| = 0.59 n$ such $E_{R \setminus S} \left(N_{R \setminus S}^t(x), B \right) = \varnothing$. It follows by our assumptions on $R$ and $S$ that a.a.s. $\textrm{deg}_G(v, B) \leq 0.51 d$ holds for every $v \in N_{R \setminus S}^t(x)$. However, by Lemma~\ref{lem::halfDegree}, this occurs with probability $o(1)$. 

Hence, given any two vertices $x, y \in V(G) \setminus S$, the above argument implies that a.a.s. $\left|N_{R \setminus S}^{t+1}(x) \right| \geq 2 n/5$ and $\left|N_{R \setminus S}^{t+1}(y) \right| \geq 2 n/5$. If $N_{R \setminus S}^{t+1}(x) \cap N_{R \setminus S}^{t+1}(y) \neq \varnothing$, then there is a path of length at most $2t+2 \leq 3 \frac{\log n}{\log (cd)}$ between $x$ and $y$ in $R \setminus S$. Otherwise, it follows by Lemma~\ref{lem::robustEdgeBetweenSets} that a.a.s. there is an edge of $R$ between $N_{R \setminus S}^{t+1}(x)$ and $N_{R \setminus S}^{t+1}(y)$, yielding a path of length at most $2t+3 \leq 3 \frac{\log n}{\log (cd)}$ between $x$ and $y$ in $R \setminus S$.
\end{proof}

The following result asserts that a relatively short path in $\mathcal{G}_{n,d}$ typically only covers a small part of the neighbourhood of every vertex outside this path.
\begin{lemma} \label{lem::fewNeighboursinPath}
Let $G \sim \mathcal{G}_{n,d}$ for some integer $d \geq 3$. Then, a.a.s. for every path $P$ in $G$ on at most $\log n/ \log d$ vertices, and every vertex $u \in V(G) \setminus V(P)$ it holds that $\emph{deg}_G(u, V(P)) \leq 3$.
\end{lemma}

Our proof of Lemma~\ref{lem::fewNeighboursinPath} makes use of the following known result regarding the edge distribution of random regular graphs; it is a rephrased version of a special case of Lemma 4.7 in~\cite{KLM}. In order to state it, we use $\mathcal{G}^*_{n,d}$ to denote the uniform probability space of random $d$-regular multigraphs with vertex-set $[n]$.
\begin{theorem} [\cite{KLM}] \label{th::GndLikeGnp}
Let $d \geq 3$ be an integer and let $E_0$ be a set of $m \leq d n/4$ pairs of elements of $[n]$. Then
$$
\mathbb{P}(E_0 \subseteq E(\mathcal{G}^*_{n,d})) \leq \left(\frac{2d}{n} \right)^m .
$$
\end{theorem}

\begin{proof} [Proof of Lemma~\ref{lem::fewNeighboursinPath}]
It is well-known (see, e.g., Theorem 9.9 in~\cite{JLR}) that if a property holds a.a.s. in $\mathcal{G}^*_{n,d}$, then it also holds a.a.s. in $\mathcal{G}_{n,d}$. Hence, it suffices to prove the corresponding claim for $\mathcal{G}^*_{n,d}$.

Let $G^* \sim \mathcal{G}^*_{n,d}$. Let $v_1, \ldots, v_k$ be a sequence of $k \leq \log n/ \log d$ vertices of $G^*$ and let $u \in V(G^*) \setminus \{v_1, \ldots, v_k\}$ be an arbitrary vertex. Let $B$ be an arbitrary subset of $\{v_1, \ldots, v_k\}$ of size $4$. Let $E_0 = \{\{v_i, v_{i+1}\} : 1 \leq i \leq k-1\} \cup \{\{u,z\} : z \in B\}$. Let $m = |E_0|$ and note that $m = k+3$. By Theorem~\ref{th::GndLikeGnp} it holds that 
$$
\mathbb{P}(E_0 \subseteq E(G^*)) \leq \left(\frac{2d}{n} \right)^m = \left(\frac{2d}{n} \right)^{k+3}.
$$
A union bound over all relevant values of $k$, and all choices of the sequence $v_1, \ldots, v_k$, the set $B$ of size $4$, and the vertex $u \in V(G^*) \setminus \{v_1, \ldots, v_k\}$, implies that the probability that there exists a path $P$ in $G^*$ on at most $\log n/ \log d$ vertices, and a vertex $u \in V(G^*) \setminus V(P)$ such that $\textrm{deg}_{G^*}(u, V(P)) \geq 4$, is at most
$$
\sum_{k=4}^{\log n/ \log d} n^k \binom{k}{4} n (2 d/n)^{k+3} \leq (\log n)^5 (2 d)^{\log n/\log d + 3} n^{-2} = o(1).
$$
\end{proof}

We can easily extend the assertion of Lemma~\ref{lem::fewNeighboursinPath} to vertices that are on the path as follows.
\begin{corollary} \label{cor::fewNeighboursinPath}
Let $G \sim \mathcal{G}_{n,d}$ for some integer $d \geq 3$. Then, a.a.s. for every path $P$ in $G$ on at most $\log n/ \log d$ vertices, and every vertex $u \in V(P)$ it holds that $\emph{deg}_G(u, V(P)) \leq 8$.
\end{corollary}

\begin{proof}
Fix a path $P = (v_1, \ldots, v_k)$ in $G$, where $k \leq \log n/\log d$, and some $1 \leq i \leq k$. Suppose for convenience that $i \notin \{1,k\}$ (the remaining cases are similar and, in fact, a bit easier). Then, $P$ can be written as $P_1 v_{i-1} v_i v_{i+1} P_2$, where $P_1$ is the (possibly empty) path $(v_1, \ldots, v_{i-2})$ and $P_2$ is the (possibly empty) path $(v_{i+2}, \ldots, v_k)$. It follows by Lemma~\ref{lem::fewNeighboursinPath} that a.a.s. $\textrm{deg}_G(v_i, V(P_1)) \leq 3$ and $\textrm{deg}_G(v_i, V(P_2)) \leq 3$. Combined with its two neighbours in $P$, namely $v_{i-1}$ and $v_{i+1}$, this yields the required bound $\textrm{deg}_G(v_i, V(P)) \leq 8$.
\end{proof}

Using Lemmas~\ref{lem::DiameterR} and~\ref{lem::fewNeighboursinPath}, and Corollary~\ref{cor::fewNeighboursinPath}, we prove a variant of Lemma~\ref{lem::S2CNP} that is suitable for random regular graphs of constant degree.
\begin{lemma} \label{lem::OurCycleLemma}
Let $G \sim \mathcal{G}_{n,d}$, where $d$ is a sufficiently large constant. Then a.a.s. the following holds. Let $R \neq G$ be a subgraph of $G$ such that $R \neq G[A, B]$ for every partition $V(G) = A \cup B$, and $\delta(R) \geq d/2$. Then, there exists an absolute constant $K$, and a cycle $C \subseteq G$ for which all of the following properties hold.
\begin{enumerate}
\item [$(a)$] $|C|$ is even and $|E(C) \setminus E(R)| = 1$;

\item [$(b)$] $|C| \leq K \frac{\log n}{\log d}$;

\item [$(c)$] $\emph{deg}_G(u, V(C)) \leq K$ holds for every $u \in V(C)$;

\item [$(d)$] $\emph{deg}_G(u, V(C)) \leq K$ holds for every $u \in V(G)\setminus V(C)$.
\end{enumerate}
\end{lemma}

\begin{proof}
Suppose that $G$ satisfies the assertions of Lemmas~\ref{lem::edgeDistribution}, \ref{lem::DiameterR}, and~\ref{lem::fewNeighboursinPath}, and of Corollary~\ref{cor::fewNeighboursinPath}; note that this fails with probability $o(1)$. 

We distinguish between the following two cases.
\begin{enumerate}
\item [(1)] $R$ is bipartite. Let $A \cup B$ be a bipartition of $R$. Since $R \neq G[A, B]$ by the premise of the lemma, there exist vertices $x \in A$ and $y \in B$ such that $xy \in E(G) \setminus E(R)$. It follows by Lemma~\ref{lem::DiameterR} that there is a path in $R$ between $x$ and $y$ of length $\ell \leq \frac{K' \log n}{\log d}$, where $K'$ is an absolute constant. Since $R$ is bipartite, $\ell$ is odd. Hence, combined with the edge $xy$, this yields a cycle $C$ satisfying properties (a) and (b). Let $r = (\ell+1) \log d/\log n = \Theta(1)$, and let $K = \max \{K', 8r\}$. Partition $C$ into $r$ consecutive vertex-disjoint paths $P_1, \ldots, P_r$, each on at most $\log n/\log d$ vertices. It follows by Lemma~\ref{lem::fewNeighboursinPath} that $\textrm{deg}_G(u, V(C)) \leq 3r$ holds for every $u \in V(G) \setminus V(C)$; this implies Property (d). Similarly, it follows by Lemma~\ref{lem::fewNeighboursinPath} and by Corollary~\ref{cor::fewNeighboursinPath} that $\textrm{deg}_G(u, V(C)) \leq 8r$ holds for every $u \in V(C)$; this implies Property (c).

\item [(2)] $R$ is not bipartite. Let $C' = (x_1, \ldots, x_{2t-1}, x_1)$ be a shortest odd cycle in $R$. Note that $|C'| \leq 3K' \frac{\log n}{\log d}$, where $K'$ is the constant whose existence is ensured by Lemma~\ref{lem::DiameterR}. Indeed, suppose for a contradiction that $|C'| > 3K' \frac{\log n}{\log d}$. By Lemma~\ref{lem::DiameterR} there is a path $P$ in $R$ between $x_1$ and $x_t$ whose length is at most $K' \frac{\log n}{\log d}$. However, $P \cup C'$ contains an odd cycle which is shorter than $C'$, contrary to the assumed minimality of $|C'|$.

Since $R \neq G$ by the premise of the lemma, there exists an edge $xy \in E(G) \setminus E(R)$. We claim that there exists a vertex $u \in V(C')$ such that there exists a path $P_x$ in $(R \setminus (V(C') \cup \{y\})) \cup \{x, u\}$ between $x$ and $u$ of length at most $K' \frac{\log n}{\log d}$. Indeed, similarly to the previous case, since $|C'| \leq 3K' \frac{\log n}{\log d}$, there exists a constant $K''$ such that $\textrm{deg}_G(u, V(C') \cup \{y\}) \leq K''$ holds for every $u \in V(G)$. It then follows by Lemma~\ref{lem::DiameterR} that the required path $P_x$ exists. An analogous argument shows that there exists a path $P_y$ in $(R \setminus (V(C') \cup V(P_x))) \cup \{y, v\}$ between $y$ and $v$ of length at most $K' \frac{\log n}{\log d}$, where $v$ is some vertex of $V(C') \setminus \{u\}$.  

Combining the path $u P_x x y P_y v$ (that has precisely one edge in $E(G) \setminus E(R)$) with one of the two paths that connect $u$ and $v$ in $C'$ (all of whose edges are in $R$) yields a cycle $C$ satisfying properties (a) and (b). As in the previous case, one can then use Lemma~\ref{lem::fewNeighboursinPath} and Corollary~\ref{cor::fewNeighboursinPath} to show that properties (c) and (d) hold as well.
\end{enumerate}
\end{proof}

Before we can state the next auxiliary result, we need the following definition.
\begin{definition} \label{def::Palphand}
A graph $G = (V, E)$ is said to have property $P_{\alpha}(n', d')$ if for every $X \subseteq V$ of size $|X| \leq n'$ and every $F \subseteq E$ such that $\left|F \cap \partial_G(x) \right| \leq \alpha \cdot \emph{deg}_G(x)$ holds for every $x \in X$, we have $|N_{G \setminus F}(X)| \geq 2 d' |X|$.
\end{definition}

The following result is essentially Theorem 3.5 from~\cite{DKN}. However, that theorem has an algorithmic component which we do not need. Therefore, we rephrase (and weaken) it here to better suit our needs. 
\begin{theorem} [Theorem 3.5 in~\cite{DKN}, abridged] \label{th::PathsInExpanders}
Let $G$ be a graph which satisfies the property $P_{\alpha}(n', d')$ for some $3 \leq d' < n'$. Suppose further that $e_G(A,B) > 0$ holds for any two disjoint sets $A, B \subseteq V(G)$ of sizes $|A|, |B| \geq n' (d' - 1)/16$. Let $S \subseteq V(G)$ be a set for which $|N_G(x) \cap S| \leq \beta \cdot \emph{deg}_G(x)$ holds for any vertex $x \in V(G)$. Let $a_1, \ldots, a_t, b_1, \ldots, b_t$ be $2t$ vertices in $S$, where $t \leq \frac{d' n' \log d'}{15 \log n'}$. If $\beta < 2 \alpha - 1$, then $G$ admits pairwise vertex-disjoint paths $P_1, \ldots, P_t$ such that for every $1 \leq i \leq t$, the endpoints of $P_i$ are $a_i$ and $b_i$.  
\end{theorem}

We are now in a position to prove Theorem~\ref{th::CycleSpaceGnd}. 

\begin{proof} [Proof of Theorem~\ref{th::CycleSpaceGnd}]
Let $G \sim \mathcal{G}_{n,d}$, where $d$ is a sufficiently large constant. Suppose that $G$ satisfies the assertions of Theorem~\ref{th::HamGnd} and of Lemma~\ref{lem::edgeDistribution}; note that this fails with probability $o(1)$. 

Suppose for a contradiction that $\mathcal{C}_n(G) \neq \mathcal{C}(G)$. We follow the recipe that was presented in Section~\ref{sec::prelim}. That is, we need to handle steps (S1), (S2), and (S3). Our assumption that $\mathcal{C}_n(G) \neq \mathcal{C}(G)$ will then lead to the contradiction appearing in (S5).  

Let $R$ be as in the premise of Lemma~\ref{lem::subgraphR}; combined with Theorem~\ref{th::HamGnd}, this takes care of (S1). Suppose that $G$  satisfies the assertions of Lemmas~\ref{lem::DiameterR} and~\ref{lem::OurCycleLemma} (with respect to $R$); note that this fails with probability $o(1)$. 

Next, we take care of (S2). Starting with (S2a), it follows by Lemma~\ref{lem::OurCycleLemma} that $G$ contains a cycle $C$ satisfying properties (a), (b), (c), and (d). In particular, $C = (v_1, \ldots, v_{2k})$ is an even cycle having an odd number of edges in $R$, and its length is at most $K \frac{\log n}{\log d}$ for some absolute constant $K$.

Prior to handling (S2b) and thinking ahead to Step (S3), we partition $Y := V(G) \setminus V(C)$ into two sets, each containing many of the neighbours of every vertex. By properties (b), (c), and (d) from Lemma~\ref{lem::OurCycleLemma} we may apply Lemma~\ref{lem::LLLsplit} to $G$ and $Y$ to obtain a set $A \subseteq Y$ of size $n/2$ and a set $B := V(G) \setminus A$ such that $\textrm{deg}_G(u, A) \geq d/9$ and $\textrm{deg}_G(u, B) \geq d/9$ hold for every $u \in V(G)$.

Returning to (S2b), let $G_1$ be the graph obtained from $G[B]$ by deleting $v_1$ and $v_{k+1}$ and all the edges (but none of the other vertices) of $C$; note that $\delta(G_1) \geq d/9 - 4 \geq d/10$. It thus follows by Lemma~\ref{lem::subgraphExpander} that $G_1$ satisfies the property $P_{2/3}(n/d, c d/2)$, where $c := c(2/3, 1/10)$ is the constant whose existence is ensured by Lemma~\ref{lem::subgraphExpander}. Moreover, by Lemma~\ref{lem::edgeDistribution}(a), there is an edge of $G_1$ between any two disjoint subsets of $V(G_1)$, each of size at least $c n/40$. Setting $S = \{v_2, \ldots, v_k, v_{k+2}, \ldots, v_{2k}\}$, observe that $|N_{G_1}(x) \cap S| \leq d/100$ holds for any vertex $x \in V(G_1)$ by properties (c) and (d) from Lemma~\ref{lem::OurCycleLemma} and by our assumption that $d$ is sufficiently large. It thus follows by Theorem~\ref{th::PathsInExpanders} that $G_1$ admits pairwise vertex-disjoint paths $P_2, \ldots, P_k$ such that, for every $2 \leq i \leq k$, the endpoints of $P_i$ are $v_i$ and $v_{2k-i+2}$. 
 
Finally, we establish (S3). Let $W = (\{v_1, \ldots, v_{2k}\} \cup V(P_2) \cup \ldots \cup V(P_k)) \setminus \{v_1, v_{k+1}\}$, and let $G_2 = G[V(G) \setminus W]$. Since $A \subseteq V(G) \setminus W$, it follows that $\delta(G_2) \geq d/9$. Therefore, by Lemma~\ref{lem::edgeDistribution}(a), by Lemma~\ref{lem::subgraphExpander}, and by our assumption that $d$ is sufficiently large, we have that $G_2$ is a $c'$-expander, where $c'$ is a sufficiently large constant, as per Theorem~\ref{th::HamConCexpander}  (note that the expansion of large sets, which are not covered by Lemma~\ref{lem::subgraphExpander}, is ensured by Lemma~\ref{lem::edgeDistribution}(a)). It then follows by Theorem~\ref{th::HamConCexpander} that $G_2$ admits a Hamilton path whose endpoints are $v_1$ and $v_{k+1}$.
\end{proof}

We end this section with a brief sketch of the proof of Theorem~\ref{th::evenn}. Let $G \sim \mathcal{G}_{n,d}$, where $n$ is even and $d \geq d_0$. For every $v \in V(G)$ let $G_v = G \setminus \{v\}$. The proof of Theorem~\ref{th::CycleSpaceGnd} carries over mutatis mutandis to show that a.a.s. $\mathcal{C}_{n-1}(G_v) = \mathcal{C}(G_v)$ holds for every $v \in V(G)$; this in turn implies that $\mathcal{C}_k(G) \subseteq \mathcal{C}_{n-1}(G)$ for every $3 \leq k \leq n-1$. Indeed, since $G$ a.a.s. satisfies the assertion of Lemma~\ref{lem::edgeDistribution} and since the removal of any single vertex from $G$ has almost no effect on the edge distribution of $G$, it follows that a.a.s. $G_v$ satisfies the assertion of Lemma~\ref{lem::edgeDistribution} for every $v \in V(G)$ (with a slightly worse error term, but the non-optimal constant 3 appearing in Lemma~\ref{lem::edgeDistribution} already accounts for this additional error). Given this fact, one may replace Theorem~\ref{th::HamGnd} with Theorem~\ref{th::HamConCexpander} (with a much larger, yet still constant, value of $d$). The remaining parts of the proof are essentially the same.

It remains to prove that a.a.s. $\mathcal{C}_n(G) \subseteq \mathcal{C}_{n-1}(G)$. Let $C = (u_1, \ldots, u_n, u_1)$ be some Hamilton cycle of $G$. Let $u_i u_j \in E(G) \setminus E(C)$ be some chord of $C$ (such an edge exists since $d \geq 3$). Let $C_1 = (u_i, \ldots, u_j, u_i)$ and let $C_2 = (u_j, u_{j+1} \ldots, u_i, u_j)$, that is, $C_1$ and $C_2$ are two non-Hamiltonian cycles in $C \cup \{u_i u_j\}$. Since $\max \{|C_1|, |C_2|\} \leq n-1$, it follows that $C_1, C_2 \in \mathcal{C}_{n-1}(G)$. Hence, $C = C_1 + C_2 \in \mathcal{C}_{n-1}(G)$ as required.

\section{Concluding remarks and open problems} \label{sec::openprob}

We have proved (see Theorem~\ref{th::CycleSpaceGnd}) that there exists a constant $d_0$ such that if $G \sim \mathcal{G}_{n,d}$ for odd $n$ and even $d \geq d_0$, then a.a.s. $\mathcal{C}_n(G) = \mathcal{C}(G)$. It would be interesting to determine how small the value of $d_0$ can be. In particular, one cannot help but wonder whether $d_0 = 4$ is sufficient (or, possibly, $d_0 = 3$ is sufficient in the setting of Theorem~\ref{th::evenn}).

An interesting extension of Theorem~\ref{th::CycleSpaceGnd} would be to prove an analogous result for $(n, d, \lambda)$-graphs. As noted in the introduction, this was previously considered in~\cite{CNP} (see Corollary 1.5 there), where it has been proved that if $G$ is an $(n, d, \lambda)$-graph, $n$ is odd, $d \geq C \log n$, and $\lambda \leq \varepsilon d/\log n$, where $C$ and $\varepsilon$ are appropriate constants, then $\mathcal{C}_n(G) = \mathcal{C}(G)$. A straightforward adaptation of our proof of Theorem~\ref{th::CycleSpaceGnd} yields the following improvement.

\begin{theorem} \label{th::ndlambda}
There exist positive constants $C$ and $\varepsilon$ such that the following holds. Let $G$ be an $(n, d, \lambda)$-graph, where $n$ is odd, $d \geq \frac{C \log n}{\log (d/\lambda)}$ is even, and $\lambda \leq \varepsilon d$. Then, $\mathcal{C}_n(G) = \mathcal{C}(G)$. 
\end{theorem}

Indeed, Theorem 1.5 in~\cite{DMMPS} may replace Theorem~\ref{th::HamGnd} (with a substantially larger but still constant $d$). Lemma~\ref{lem::LLLsplit} and Theorem~\ref{th::PathsInExpanders} are unaffected. Lemmas~\ref{lem::edgeDistribution}, \ref{lem::halfDegree}, and~\ref{lem::robustEdgeBetweenSets} remain essentially the same. In Lemma~\ref{lem::subgraphExpander} we will get an expansion by a factor of $c d^2/ \lambda^2$ (instead of the current $c d$). This (in general weaker) expansion will lead to paths of length at most $\frac{c' \log n}{\log (d/\lambda)}$ in Lemma~\ref{lem::DiameterR}. Skipping Lemma~\ref{lem::fewNeighboursinPath}, Theorem~\ref{th::GndLikeGnp}, and Corollary~\ref{cor::fewNeighboursinPath}, in which we needed true randomness, we obtain analogues of properties (c) and (d) in Lemma~\ref{lem::OurCycleLemma} simply because the cycle we construct will have length less than, say, $d/100$. The remainder of the proof is essentially the same. 

It would be interesting to determine whether, as in the case of random regular graphs, an analogue of Theorem~\ref{th::ndlambda} could hold for constant $d$ and $\lambda \leq \varepsilon d$.

Theorem~\ref{th::perturbedHam}, and thus also Theorem~\ref{th::CycleSpacePerturbed}, are essentially best possible. Nevertheless, in~\cite{BFM} Bohman, Frieze, and Martin proved that adding a mild upper bound on the independence number of the graph being perturbed drastically decreases the amount of random perturbation needed to ensure Hamiltonicity (see also~\cite{AHK} for various extensions and related results). More formally, if $H$ is an $n$-vertex graph with minimum degree $\delta(H) > \delta n$ and independence number $\alpha(H) \leq \delta^2 n/2$, and $G \sim \mathbb{G}(n,p)$, where $p := p(n) = \omega \left(n^{-2} \right)$, then $H \cup G$ is a.a.s. Hamiltonian. It would be interesting to determine whether the same (or similar) conditions are sufficient to ensure that a.a.s. $\mathcal{C}_n(H \cup G) = \mathcal{C}(H \cup G)$.


\begin{thebibliography}{99}

\bibitem{AHK}
E. Aigner-Horev, D. Hefetz, and M. Krivelevich, Cycle lengths in randomly perturbed graphs, \emph{Random Structures and Algorithms} 63(4) (2023), 867--884.

\bibitem{AS}
N. Alon and J. H. Spencer, \textbf{The Probabilistic Method}, Wiley-Interscience Series in Discrete Mathematics
and Optimization, John Wiley and Sons, fourth edition, 2015.

\bibitem{ALW}
B. Alspach, S. C. Locke, and D. Witte, The Hamilton spaces of Cayley graphs on abelian groups, \emph{Discrete Mathematics} 82 (2) (1990), 113--126.

\bibitem{BK}
J. D. Baron and J. Kahn, On the cycle space of a random graph, \emph{Random Structures and Algorithms} 54(1) (2019), 39--68.

\bibitem{BFM}
T. Bohman, A. Frieze, and R. Martin, How many random edges make a dense graph Hamiltonian?, \emph{Random Structures and Algorithms} 22 (2003), 33--42.

\bibitem{BL}
J. A. Bondy and L. Lov\'asz, Cycles through specified vertices of a graph, \emph{Combinatorica} 1 (1981), 117--140.

\bibitem{CNP}
M. Christoph, R. Nenadov, and K. Petrova, The Hamilton space of pseudorandom graphs, arXiv preprint arXiv:2402.01447, 2024.

\bibitem{DHJ}
B. DeMarco, A. Hamm, and J. Kahn, On the triangle space of a random graph, \emph{Journal of Combinatorics} 4(2) (2013), 229--249.

\bibitem{Dirac}
G. A. Dirac, Some theorems on abstract graphs, \emph{Proceedings of the London Mathematical Society} 3(1) (1952), 69--81.

\bibitem{DKN}
N. Dragani\'c, M. Krivelevich, and R. Nenadov, Rolling backwards can move you forward: on embedding problems in sparse expanders, \emph{Transactions of the American Mathematical Society} 375 (7) (2022), 5195--5216.

\bibitem{DMMPS}
N. Dragani\'c, R. Montgomery, D. Munha Correia, A. Pokrovskiy, and B. Sudakov, Hamiltonicity of expanders: optimal bounds and applications, arXiv preprint arXiv:2402.06603v2, 2024.

\bibitem{Friedman}
J. Friedman, A proof of Alon’s second eigenvalue conjecture and related problems, \emph{Memoirs of the American Mathematical Society} 195 (910), (2008).

\bibitem{HMMMP}
M. Hahn-Klimroth, G. S. Maesaka, Y. Mogge, S. Mohr, and O. Parczyk, Random perturbation
of sparse graphs, \emph{Electronic Journal of Combinatorics} 28 (2021), Paper 2.26.

\bibitem{Hartman}
I. B.-A. Hartman, Long cycles generate the cycle space of a graph, \emph{European Journal of Combinatorics} 4 (1983), 237--246.

\bibitem{HK}
D. Hefetz and M. Krivelevich, The Hamilton cycle space of random graphs, arXiv preprint arXiv:2506.19731v1, 2025.

\bibitem{HKT}
D. Hefetz, M. Krivelevich and T. Szab\'o, Sharp threshold for the appearance of certain spanning trees in random graphs,
\emph{Random Structures and Algorithms} 41 (2012), 391--412.

\bibitem{Heinig1}
P. Heinig, On prisms, M\"obius ladders and the cycle space of dense graphs, \emph{European Journal of
Combinatorics} 36 (2014), 503--530.

\bibitem{Heinig2}
P. Heinig, When Hamilton circuits generate the cycle space of a random graph, arXiv preprint arXiv:1303.0026, 2013.

\bibitem{HY}
X. Hou and Z. Yin, Dirac-type condition for Hamilton-generated graphs, arXiv preprint arXiv:2503.15950v1, 2025.

\bibitem{JLR}
S. Janson, T. \L uczak, and A. Ruci\'nski, \textbf{Random graphs}, Wiley-Interscience Series in Discrete Mathematics and Optimization, Wiley-Interscience, New York, 2000.

\bibitem{KLM}
M. Krivelevich, A. Lew, and P. Michaeli, Rigid partitions: from high connectivity to random graphs, arXiv preprint arXiv:2311.14451v2, 2023.

\bibitem{Locke1}
S. C. Locke, A basis for the cycle space of a 2-connected graph, \emph{European Journal of Combinatorics} 6 (1985), 253--256.

\bibitem{Locke2}
S. C. Locke, A basis for the cycle space of a 3-connected graph, \emph{Annals of Discrete Mathematics} 27 (1985), 381--397.

\bibitem{RW1}
R. W. Robinson and N. C. Wormald, Almost all cubic graphs are Hamiltonian, \emph{Random Structures and Algorithms} 3 (1992), 117--126.

\bibitem{RW2}
R. W. Robinson and N. C. Wormald, Almost all regular graphs are Hamiltonian, \emph{Random Structures and Algorithms} 5 (1994), 363--374.

\end{thebibliography}
\end{document}